\documentclass[a4paper,11pt]{amsart}

\usepackage{amsmath}
\usepackage{amssymb}
\usepackage{amstext}
\usepackage{amsbsy}
\usepackage{amsopn}
\usepackage{amsthm}
\usepackage{amscd}
\usepackage{amsxtra}
\usepackage{amsfonts}

\usepackage{ifthen}
\usepackage{xspace}
\usepackage{fancyheadings}
\usepackage{layout}
\usepackage{hyperref}

\newcommand{\preprintserver}[2]{\href{http://arXiv.org/abs/math/#2}{#1/#2}}

\input xy
\xyoption{matrix}
\xyoption{arrow}
\xyoption{curve}
\newdir{ (}{{}*!/-5pt/@^{(}}

\newboolean{xlabels}
\newcommand{\xlabel}[1]{
                        \label{#1}
                        \ifthenelse{\boolean{xlabels}}
                                   {\marginpar{#1}}
                                   {}
                       }

\newcommand{\ZZ}{\mathbb{Z}}

\newcommand{\AZ}{\mathbb{A}}

\newcommand{\CC}{\mathbb{C}}
\newcommand{\RR}{\mathbb{R}}
\newcommand{\QQ}{\mathbb{Q}}
\newcommand{\NN}{\mathbb{N}}
\newcommand{\PP}{\mathbb{P}}

\newcommand{\FF}{\mathbb{F}}

\newcommand{\suchthat}{\, | \,}

\newboolean{probleme}
\newcommand{\problem}[1]
           {\ifthenelse{\boolean{probleme}}
                       {{\bf(PROBLEM: #1)\bf}}
                       {}
           }
\newboolean{zukuenftiges}
\newcommand{\zukunft}[1]
           {\ifthenelse{\boolean{zukuenftiges}}
                       {{\bf(AUSBAUM\"OGLICHKEIT: #1)\bf}}
                       {}
           }
\newboolean{extras}
\newcommand{\extra}[1]
           {\ifthenelse{\boolean{extras}}
                       {{\bf EXTRA #1 EXTRA\bf}}
                       {}
           }

\newboolean{ignore}
\newcommand{\ignore}[1]
           {\ifthenelse{\boolean{ignore}}
                       {{\bf IGNORE #1 IGNORE\bf}}
                       {}
           }

\DeclareMathOperator{\codim}{codim}

\DeclareMathOperator{\SO}{SO}

\DeclareMathOperator{\Char}{char}

\theoremstyle{plain}
\begingroup

\newtheorem{thm}{Theorem}
\newtheorem{cor}[thm]{Corollary}
\newtheorem{lem}[thm]{Lemma}

\newtheorem{prop}[thm]{Proposition}
\newtheorem{conj}[thm]{Conjecture}

\numberwithin{thm}{subsection} 
\endgroup

\newtheorem*{thm*}{Theorem}
\newtheorem*{conj*}{Conjecture}
\newtheorem*{verm*}{Vermutung}

\theoremstyle{definition}
\newtheorem{defn}[thm]{Definition}
\newtheorem{rem}[thm]{Remark}
\newtheorem{example}[thm]{Example}

\newtheorem{notation}[thm]{Notation}

\newtheorem{experiment}[thm]{Experiment}
\newtheorem{alg}[thm]{Algorithm}
\newtheorem{heu}[thm]{Heuristic}
\newtheorem{cau}[thm]{Caution}

\numberwithin{equation}{section}

\newcommand{\nosubsections}{\renewcommand{\thethm}{\thesection.\arabic{thm}}
                            \setcounter{thm}{0}
                           }

\newcommand{\cref}[3]{(\ref{#1}, #2 \ref{#3})}

\date{\today}

\setboolean{probleme}{true}
\setboolean{xlabels}{false}

\usepackage{graphicx}

\newcommand{\zoladek}{\.Zo\l\c adek\,}
\newcommand{\zoladeks}{\.Zo\l\c adek's\,}

\begin{document}

\title{Experimental results for the Poincar\'e center problem}

\address{Institiut f\"ur Mathematik\\
          Universit\"at Hannover\\
          Welfengarten 1\\
          D-30167 Hannnover
         }


\urladdr{\href{http://www-ifm.math.uni-hannover.de/~bothmer}{http://www-ifm.math.uni-hannover.de/\textasciitilde bothmer}}


\author[v. Bothmer]{Hans-Christian Graf v. Bothmer}

\begin{abstract}
We apply a heuristic method based on counting points over finite fields to the Poincar\'e center problem.
We show that this method gives the correct results for homogeneous non linearities of degree $2$ 
and $3$. Also we obtain new evidence for \zoladeks conjecture about general degree $3$ non linearities.
\end{abstract}

\maketitle

\section*{Introduction}
\nosubsections

\newcommand{\Fp}{\FF_p}
\newcommand{\Anp}{\AZ^n(\Fp)}
\newcommand{\gammatilde}{\tilde{\gamma}}

In 1885 Poincar\'e asked when the differential equation
\[
y' = - \frac{x + p(x,y)}{y+q(x,y)} =: - \frac{P(x,y)}{Q(x,y}
\]
with convergent power series $p(x,y)$ and $q(x,y)$ starting with quadratic terms, 
has stable solutions in the neighborhood of the equilibrium solution
$(x,y)=(0,0)$. This means that in such a neighborhood the solutions of the
equivalent plane autonomous system
\begin{align*}
	\dot{x} &= y + q(x,y) = Q(x,y)\\
	\dot{y} &= -x - p(x,y) = -P(x,y)
\end{align*}
are closed curves around $(0,0)$.

Poincar\'e showed that one can iteratively find a formal power series
$F = x^2+y^2+f_2(x,y)+f_3(x,y)+\dots$ such that
\[
	\det \begin{pmatrix} F_x & F_y \\ P & Q \end{pmatrix} = \sum_{j=1}^\infty s_j(x^{2j+2}+y^{2j+2})
\]
with $s_j$ polynomials in the coefficients of $P$ and $Q$.
If all $s_j$ vanish, and $F$ is convergent then $F$ is a constant of motion, i.e. its gradient field
satisfies $Pdx+Qdy=0$. Since $F$ starts with $x^2+y^2$ this shows that close to the origin all integral curves are closed and the system is stable. Therefore the $s_j$'s are
called the {\sl focal values} of $Pdx+Qdy$.

Poincar\'e also showed, that if an analytic constant of motion exists, the focal values must vanish.
Later Frommer \cite{Frommer} proved that the systems above are stable if and only if all focal values vanish even without the assumption of convergence of $F$. (Frommer's proof contains a gap which can be closed \cite{vWahlGap})

Unfortunately it is in general impossible to check this condition for a given differential equation because
there are infinitely many focal values. In the case where $P$ and $Q$ are polynomials of degree
at most $d$, the $s_j$ are polynomials in finitely many unknowns. Hilbert's Basis Theorem then implies
that the ideal $I_\infty = (s_1,s_2,\dots)$ is finitely generated, i.e there exists an integer $m := m(d)$
such that
\[
		s_1 = s_2 = \dots = s_{m(d)} = 0  \implies s_j = 0 \quad\forall j.
\]
This shows that a finite criterion for stableness exists, but due to the indirect proof of
Hilbert's Basis Theorem no value for $m(d)$ is obtained. In fact even today only $m(2)=3$ is known. 
\zoladek \cite{ZoladekEleven} and Christopher \cite{ChristopherEleven} showed that $m(3) \ge 11$.
Our experimental data at this moment recovers $m(3) \ge 10$.

The proof for $m(2)=3$ is conceptually simple: Compute the first $3$ focal values as polynomials
in the coefficients of $P$ and $Q$ under the assumption $\deg(P)=\deg(Q)=2$. The $3$ polynomials cut out an algebraic variety in the space of all differential equations
of degree $2$. Then decompose, by hand or by computer, this variety into its irreducible components.
For each component prove that all its differential equations have a constant of motion. By works of
Dulac \cite{Dulac} and Schlomiuk \cite{SchlomiukTransactions} we know that in this case integrating factors of Darboux type, i.e of the form $\prod_{i=1}^{k} F_i^{\alpha_i}$ with $F_i$ polynomial always suffice to find constants of motion.

For $d=3$ this approach is not feasible because the polynomials $s_j$
are very large. They involve $14$ variables and are of weighted degree $2j$. 
For example $s_5$ has already $5348$ terms and takes about $1.5$ hours on a
Powerbook G4 to calculate. The polynomials $s_j$, $j\ge 6$ can at the moment not be determined by
computer algebra systems. Even if we would somehow obtain these polynomials, it is
extremely difficult to decompose the resulting variety into irreducible components. Even $I_5 = (s_1,\dots,s_5)$ can not be decomposed by current systems. So for $d=3$ only partial results are known. 

One new feature of the center problem in degree $d=3$ is that a new type of differential equations with
a center appear, namely the rationally reversible centers.
\zoladek has classified all rational reversible cubic systems which are not of Darboux type in \cite{ZoladekRational} and \cite{ZoladekCorrection}. He also conjectured that all cubic systems with stable solutions near the origin are either of Darboux type or rationally reversible. This conjecture has been verified on several linear subspaces of the space of all differential forms, for example in \cite{ZoladekRomanovskii} and \cite{ChristopherSubspace}. In this paper we provide additional statistical evidence for \zoladeks conjecture in the whole space of Poincar\'e differential equations up to codimension $7$.

Our main tool is a statistical method of Schreyer \cite{irred} to estimate the number
of components of the locus $X_i$ where the first $i$ focal values vanish. 
The basic idea is to reduce the
equations $s_k$ modulo a prime number $p$ and count the number
of $\Fp$-rational points of $X_i$. By the Weil Conjectures \cite{weilconjectures}, which were proved by Delinge \cite{deligneproof}, we know that
the fraction of points 
\[
	\gamma_p(X_i) := \frac{\text{number of $\Fp$ rational points on $X_i$}}
					{p^{14}}
\]
is equal to 
\[
	r\Bigl( \frac{1}{p} \Bigr)^c + \text{higher order terms}
\]
where $c$ is the maximal codimension
of a component of $X$ and $r$ is the number of components with this codimension. 

To evaluate the $s_k$ at all possible points is not feasible, but by using a large number of
random points one obtains an approximate value of $\gamma_p(X)$ that can be used to 
estimate $r$ and $c$. 

In an appendix to this paper which is joint work with Martin Cremer we show how one can use
Frommer's algorithm to evaluate the $s_k$'s in a given point over a finite field $\Fp$ without knowing the explicit polynomial equations for $s_k$. We describe a C++ implementation following 
 \cite{Frommer}, \cite{Moritzen} and \cite{Hoehn} that can evaluate 3.000.000 points per hour on a 450 MHz machine. 

As a slight improvement to Schreyer's method we also look at the tangent space of $X_i$ in
the random points. In the appendix we also show that these tangent spaces can be calculated using Frommer's algorithm. By using the inequality
\[
	\codim (T_{X,x}) \le  \codim X'
\]
where $X'$ is a component of $X$ passing through $x$ we can eliminate all points on components
of codimension at most $c$ when estimating the number of components in codimension $c+1$. 
This allows us to estimate also the number of components of larger codimension. 

The method of investigating many random points has minimal memory requirements and can be parallelized with almost no overhead. It also has the advantage of giving continuous intermediate results while running and can therefore be stopped and restarted with no loss.

This characteristic $p$ method differs from the one applied by Fronville \cite{Fronville}. She reduces explicit focal polynomials modulo $p$ to make  factorisation and Gr\"obner basis computations faster. For general cubic centers the explicit polynomials of high degree are unknown so Gr\"obner basis and factorisation methods can not be applied. Even if the polynomials where known, they are to large to be factored with current methods, even over a finite field.

We start by applying our method to homogeneous differential equations of degree $2$ and $3$ where we recover the results known from literature. For general degree $3$ differential equations we obtain

\begin{conj*}
Let $X_{\infty} \subset \AZ^{14}$ be the algebraic set of inhomogeneous 
degree $3$ Poincar\'e differential equations
over $\CC$ where all focal values vanish. Then $X_\infty$ has

\begin{enumerate}
\item $1$ component of codimension $5$, 
\item $2$ components of codimension $6$,
\item $4$ components of codimension $7$,
\item at least one component of codimension $8$,
\item unknown numbers of components of codimension $\ge 9$.
\end{enumerate}

\end{conj*}

In Section \ref{sAlgebraic} we investigate the known families of Darboux type Poincar\'e differential equations. Using a theorem of Christopher we exhibit Darboux type families whose closures
form  components of $X_\infty$ with codimensions $5, 6, 7, 7$ and $7$. These are those with algebraic integral curves of degrees $(4), (3,1), (2,2), (2,1,1)$ and $(1,1,1,1)$ respectively.

In Section \ref{sReversible} we show that there exist two families of rational reversible systems
whose closures form components of $X_\infty$ of codimension $6$ and $7$. Namely these are the
symmetric systems and \zoladeks family $CR_{11}$. Together with these results our experiments
suggest that these components form a complete list up to codimension $7$ and that \zoladeks conjecture holds up to this codimension. 

As a further check of our experiment we look a the other families of \zoladeks list in Section \ref{sList}
and show that
all of them either lie on one of the above mentioned components or on components of codimension at least $8$. For this we find new invariant conics for \zoladeks families $CR_5, CR_7, CR_{12}$ 
and $CR_{16}$ and show that general members of these families are also of Darboux type.

To conclude we show in Section \ref{sNecessary} how our methods can be used for a heuristic estimate of the number of equations needed to define $X_{\infty}$. This leads to the correct estimates for homogeneous equations of degree $2$
and $3$. For general degree $3$ equations we do not yet have enough data to give a reliable estimate. Our experiments so far give a lower bound $m(3) \ge 10$, which doesn't quite reach the known bound $m(3)\ge 11$. To obtain heuristic evidence for a possible upper bound $m(3) \le 11$ we would need to speed up our calculations by a factor of 10000.

I would like to thank Wolf v.\,Wahl for bringing the Poincar\'e center problem to my attention and for the helpful discussions about this topic. Also a discussion with S\l awomir Cync about real and complex solution sets was very useful. Furthermore I would like to thank Colin Christopher for pointing out \zoladeks results to me.

Part of the calculations of this article were done at the Gauss Laboratory at the University of G\"ottingen. I would like to thank Yuri Tschinkel for providing this opportunity.

\section{Preliminaries}
\nosubsections

In this paper we write the differential equation
\[
y' = - \frac{P(x,y)}{Q(x,y)}
\]
as $P(x,y)dx + Q(x,y)dy = 0$. If $P$ and $Q$ are polynomials of degree at most $d$ we can homogenize
$P$ and $Q$ with respect to a third variable $z$. Conversely we obtain a differential equation
\[
y' = - \frac{P(x,y,1)}{Q(x,y,1)}
\]
from every homogeneous polynomial differential form $P(x,y,z)dx + Q(x,y,z)dy$. If not stated differently we will use the homogeneous formulation. 

\newcommand{\Vd}{V^d}
\newcommand{\Vdthree}{V^3}
\newcommand{\Pd}{V^d_{\text{Poincar\'e}}}
\newcommand{\Pdthree}{V^3_{\text{Poincar\'e}}}
\newcommand{\lineinf}{\PP^1_\infty}
\begin{notation}
Furthermore we denote by

\begin{tabular}{ll}
$K$ & an commutative ring (usually a field or $\ZZ$),\\
$\PP^2:=\PP^2_K$ & the projective plane over this field, \\
$K[x,y,z] $ & the coordinate ring of $\PP^2$,\\
$\AZ^2 \subset \PP^2$ & the affine plane where $z\not=0$, \\
$\lineinf \subset \PP^2$ & the line at infinity where $z=0$,\\
$K[x,y,z]_d$& the vector space of homogeneous degree $d$ polynomials,\\
$P,Q \in K[x,y,z]_d$ & two such polynomials,\\
$Pdx + Qdy$ & the corresponding differential form on $\PP^2$,\\
$\Vd \cong K[x,y,z]_d \oplus K[x,y,z]_d$ & the vector space of all such differential forms,\\
$\{ x^iy^jz^{d-i-j}dx, x^iy^jz^{d-i-j}dy \}_{i+j\le d} $& the monomial basis of this vector space,\\
$p_{ij}, q_{ij} \in K$ & the coordinates of $\Vd$ with respect to this basis. \\
\end{tabular}

Also we denote for any polynomial $F \in K[x,y,z]$ the (formal) partial derivatives
by $F_x$, $F_y$ and $F_z$.

\end{notation}

The differential equations considered by Poincar\'e have a special form:

\begin{defn}
Let $Pdx+Qdy \in V_d$ be a polynomial differential form of degree $d$ on $\PP^2$. A point $a \in \PP^2$ is called a {\sl critical point} of $Pdx+Qdy$, if $P(a)=Q(a)=0$. $Pdx+Qdy$ is called 
a {\sl Poincar\'e differential form}, if $a=(0:0:1)$ is a critical point, $P_x(a) = Q_y(a) = 1$ and $P_y(a) = Q_x(a)=0$. We denote the
vector space of all Poincar\'e differential forms by $\Pd$. $\Pd$ is a codimension $6$ linear subspace of $\Vd$.
\end{defn}

\begin{rem}
The Poincar\'e differential forms $Pdx + Qdy \in \Pd$ correspond to the differential equations
\[
	y' = -\frac{x + p(x,y)}{y + q(x,y)}
\]
of the Poincar\'e center problem, since the conditions above imply that $P$ and $Q$
are of the form
\[
	P = xz^{d-1} + \text{terms with fewer $z's$}
	\quad \quad \text{and} \quad \quad
	 Q = yz^{d-1} + \text{terms with fewer $z's$}
\]
\end{rem}

\begin{defn}
Let $Pdx + Qdy$ be a Poincar\'e differential form of degree $d$ over a field of characteristic $0$. One can then use Frommer's algorithm \ref{aFrommersAlgorithm} in the appendix to find a formal 
power series $F \in K[[x,y]]$ with
\[
	\det \begin{pmatrix} F_x(x,y) & F_y(x,y) \\ P(x,y,1) & Q(x,y,1) \end{pmatrix} = \sum_{j=1}^\infty s_j(P,Q)(x^{2j+2}+y^{2j+2}).
\]
In this situation $s_j(P,Q)$ is called the {\sl $j$th focal value} of $Pdx+Qdy$. By Corollary \ref{cPolynomial} in the appendix Frommer's algorithm also
implies that $s_j$ is polynomial on $\Pd$ and has rational coefficients. We call $s_j \in \QQ[p_{ij},q_{ij}]$ the {\sl $j$th focal
polynomial}.
\end{defn}

\begin{example} \xlabel{eFirstFocal}
For Poincar\'e differential equations of degree $3$ the first focal polynomial is:
\[
{\textstyle s_1 = \frac{2}{3} p_{02} q_{02}+\frac{1}{3} p_{02} p_{11}-\frac{1}{3} q_{02} q_{11}+\frac{1}{3} p_{11} p_{20}-\frac{1}{3} q_{11} q_{20}-\frac{2}{3} p_{20} q_{20}-p_{03}+\frac{1}{3} q_{12}-\frac{1}{3} p_{21}+q_{30} }
\]
The size of the focal polynomials grows very fast. With current computer algebra systems one
can not calculate $s_k$ with $k\ge6$ for degree $3$ Poincar\'e differential equations. For this reason we do not use explicit polynomials in this article.
\end{example}

\begin{defn}
For $s_k \in \QQ[p_{ij},q_{ij}]$ we denote by $\delta_k$ the smallest common denominator
of all coefficients, i.e. $\delta_k s_k \in \ZZ[p_{ij},q_{ij}]$. We call $I_i = (\delta_1 s_1, \dots \delta_i s_i)$
the $i$-th {\sl center ideal} and $I_\infty = (\delta_i s_i)_{i \in \NN}$ the {\sl total center ideal}.

The vanishing sets of these ideals are the $i$-th {\sl center variety} $X_i = V(I_i) \in \Pd$
and the {\sl total center variety} $X_\infty = V(I_\infty) \subset \Pd$. Since $I_i, I_\infty \subset \ZZ[p_{ij},q_{ij}]$ these varieties are defined over every commutative ring $K$.
\end{defn}

\begin{rem}
By Corollary \ref{cFiniteField} in the appendix on can use Frommer's algorithm over a finite field to check $\delta_js_j(P,Q)=0$  without knowing $\delta_js_j$ explicitly. 
By Remark \ref{rEpsilon} one can use Frommer's algorithm over $K[\epsilon]/(\epsilon^2)$ to determine the tangent space
to $X_j$ in a given point $Pdx+Qdy \in X_j$. For this it is again not necessary to know $\delta_js_j$ explicitly. 
\end{rem}

\section{Counting Points}
\nosubsections

In this section we explain how one can obtain heuristic information about a variety $X \subset \AZ^n$
by evaluation its defining equations at random points. For an extended discussion about this method see \cite{irred}.

\begin{defn}
Let $X \subset \Anp$ be an algebraic variety. Denote the number of $\Fp$-rational
points of $X$ by $|X(\Fp)|$.  Then
\[
	\gamma_p(X) = \frac{|X(\Fp)|}{|\Anp|}
\]
is called the {\sl fraction of $\Fp$-rational points of $X$ in $\AZ^n$}.
\end{defn}

\begin{rem}
If $X$ has $r$ irreducible components of codimension $c$ and all other irreducible components
have larger codimension then the Weil-Conjectures imply that
\[
	\gamma_p(X) = r \Bigl( \frac{1}{p} \Bigr)^c + \text{higher order terms in $\frac{1}{p}$}
\]
\end{rem}

We will estimate $\gamma_p(X)$ statistically by evaluating the equations defining $X$ in a
number of randomly chosen points.

\begin{defn}
Let $X \subset \Anp$ be an algebraic variety. For a sequence 
$S= (x_1,\dots,x_N)$ of $\Fp$-rational points in $\Anp$ we call
\[
	\gammatilde_p(X,S) = \frac{| \{ i \suchthat x_i \in X\}|}{N}
\]
the {\sl empirical fraction of  $\Fp$-rational points}.
\end{defn}

\begin{rem}
The distribution of $\gammatilde_p(X,S)$ on the set of all sequences $S$ of length $N$
is binomial with mean $\mu(\gammatilde_p(X,S)) = \gamma_p(X)$ and standard deviation
\[
	\sigma(\gammatilde_p(X,S)) = \sqrt{\frac{\gamma_p(X)(1-\gamma_p(X))}{N}}
\]
\end{rem}

This allows us to obtain an estimate of $\gamma_p(X)$ and then of $r$ and $c$ by evaluating
the equations of $X$ in many random points. More information is obtained, if we also calculate the
tangent space of $X$ in these random points: 

\begin{rem} 
Let $X' \subset X \subset \AZ^n$ be an irreducible component, 
$x \in X'$ a point and $T_{X',x}$ the
tangent space of $X'$ in $x$. Then
\[
		\codim X' \ge \codim T_{X',x}
\]
with equality for general points if $X'$ is reduced. We therefore consider only points with $\codim T_{X',x} = c$ in estimating the number of components of codimension $c$. By the inequality
above we disregard all points on components of codimension greater then $c$.
\end{rem}

These arguments lead us to

\begin{heu}
Evaluate the equations of $X$ in $N$ random points $x_i$ over $\FF_p$ and calculate the
tangent spaces $T_{X,x_i}$ in these points. Then estimate
\[
	\#\{ \text{codim $c$ components} \} \approx \frac{\#\{i \suchthat \codim T_{X,x_i} = c\}}{N}p^c
\]
with an estimated error
\[
	\sigma \frac{\sqrt{\#\{i \suchthat \codim T_{X,x_i} = c\}}}{N}p^c.
\]
In this paper we have used $\sigma=2$.
\end{heu}

\begin{cau}
Let $X^c$ be the subvariety of $X$ whose points have a tangent space of codimension $c$. Then above heuristic means that statistically the hypothesis $\gamma_p(X^c) = r (1/p)^c$ can not be rejected 
with confidence of more than $4.6 \%$. Algebraically this proves nothing, but gives a way to arrive at a reasonable conjecture about $X$. 
\end{cau}

\begin{rem}
One might propose to estimate the higher coefficients in the power series of the Weil formulas using the same methods as above. Notice that in this case the error scales with $(1/p)^{c_{min}}$ where
$c_{min}$ is the codimension of the largest component. Using the infinitesimal information we have the
better scaling $(1/p)^c$. In our Experiment \ref{eGen3} the second method is
$23^2 \approx 500$ times faster when estimating the number of codimension $7$ components.
\end{rem}

\section{Experiments}
\nosubsections

To show that our heuristic gives useful results we applied it to the Poincar\'e center problem with
$d=2$, $d=3$ with $p_{ij}=q_{ij}=0$ for $i+j=2$ (homogeneous case) and $d=3$ without restrictions.
For the first two cases the decomposition of $X_\infty$ into irreducible components is known and agrees with the estimates of our heuristic. In the general $d=3$ case our computer power is enough to
estimate the number of components with codimension at most $7$.

\begin{figure}
\begin{center}
{\bf \Large Homogeneous Poincar\'e differential equations}
\end{center}
\vspace{5mm}

\includegraphics[width=14cm,clip=true]{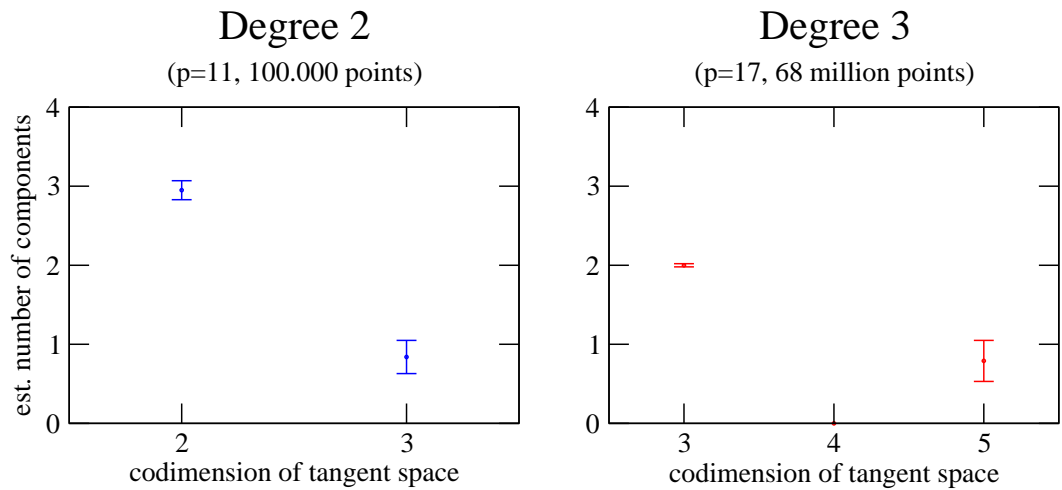}
\caption{Our heuristic gives the correct results for homogeneous differential equations. \label{fHom}}
\end{figure}

\begin{experiment} \xlabel{eHom2}
Over $\FF_{11}$ we calculated the first $4$ focal values of $100,000$ 
differential equations  with homogeneous nonlinearities of degree $2$ and random coefficients using Frommer's Algorithm as described in the appendix. If all focal values vanished, we also
calculated the codimension of the tangent space of $X_{7}$ in this point. The
results are collected in Figure \ref{fHom}. Since we found not enough points with a tangent space
of codimension $0$ or $1$ to justify a component of this codimension we show them together with the codimension $2$ points.

Notice that our results agree with the well known fact
that $X_{\infty}$ has $4$ components in this case: The Hamiltonian differential equations (codim $2$), the symmetric differential equations (codim $2$), the component of differential equations with three integral lines (codim $2$) and the component with an integral conic and and integral cubic in special position (codim $3$) \cite{SchlomiukTransactions}.
\end{experiment}

\begin{experiment} \xlabel{eHom3}
Over $\FF_{17}$ we calculated the first $7$ focal values of $68$ million
differential equations with homogeneous nonlinearities of degree $3$ and random coefficients. If all focal values vanished, we also
calculated the codimension of the tangent space of $X_{7}$ in this point. The
results are collected in Figure \ref{fHom}. Again we have shown the points with tangent spaces of codimension less then $3$ together with those of codimension $3$.

Malkin \cite{MalkinDeg3} has first given center conditions in this case. More explicitly Lunkevich and Sibirski{\u\i} have shown in \cite{LSHom3} that $I_{\infty}$ decomposes into two ideals which they call condition I and II. Condition I defines a codimension $3$ variety which on closer inspection decomposes further into the variety of Hamiltonian differential equations and the variety of symmetric differential equations. Condition II defines an irreducible variety of codimension $5$. This agrees with our experiment.\end{experiment}

\begin{figure}
\begin{center}
{\bf \Large General Poincar\'e differential equations of degree $3$}
\end{center}
\vspace{5mm}

\includegraphics[width=14cm,clip=true]{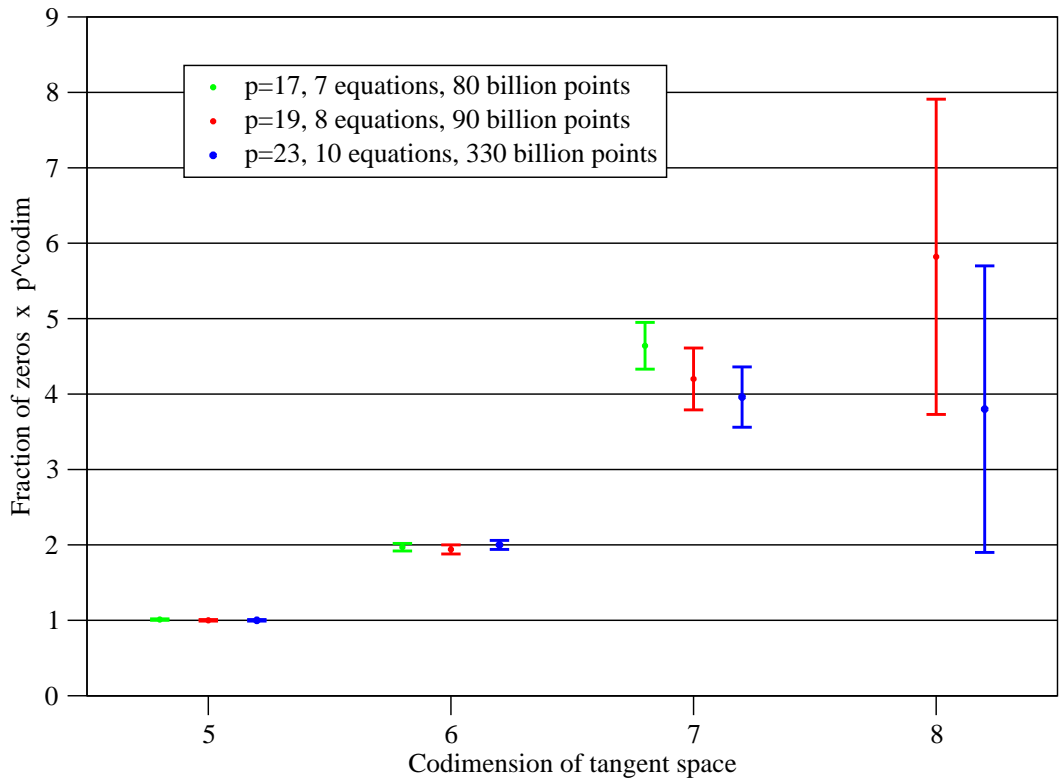}
\caption{The higher result for $p=17$ in codimension $7$ is expected, since in this case
we only consider $X_7$. The $8$th focal value could be nonzero on a general point of one component.
Similarly we expect a higher result for $p=19$ in codimension $8$. In codimension $8$ we
have only $2$ measurements, since one can not evaluate $\delta_8s_8$ over $\FF_{17}$ with Frommer\'s algorithm. \label{fGen3}}
\end{figure}

\begin{experiment} \xlabel{eGen3}
Over $\FF_p$ we calculated the first $(p-3)/2$ focal values of several billion
differential equations with general nonlinearities of degree at most $3$ and random coefficients. If all focal values vanished, we also
calculated the codimension of the tangent space of $X_{(p-3)/2}$ in this point. The
results are collected in Figure \ref{fGen3}.

\end{experiment}

This experiment is consistent with the following conjecture

\begin{conj}
Let $X_{\infty} \subset \AZ^{14}$ be the algebraic set of inhomogeneous 
degree $3$ Poincar\'e differential equations
over $\CC$ where all focal values vanish. Then $X_\infty$ has

\begin{enumerate}
\item $1$ component of codimension $5$, 
\item $2$ components of codimension $6$,
\item $4$ components of codimension $7$,
\item at least one component of codimension $8$,
\item unknown numbers of components of codimension $\ge 9$.
\end{enumerate}

\end{conj}

In the following sections we will identify known families of differential forms with a center whose closure
is a component  of $X_{\infty}$ with codimension at most $7$. The number of such families agrees with our experimental results. This provides new evidence for \zoladeks conjecture up to codimensin $7$.

\section{Darboux Centers} \xlabel{sAlgebraic}
\nosubsections

Let $Pdx+Qdy$ be a differential form of degree $d$ with algebraic integral
curves whose degrees sum up to $d+1$. If these integral curves are in general position it is well known that $Pdx+Qdy$ admits an integrating factor \cite{enoughCurves}. If $Pdx+Qdy$ is
a Poincar\'e differential form this means that $Pdx+Qdy$ has closed integral curves near the origin. 

In this section we consider for each partition $\mu$ of $d+1$ the variety $X_\mu$ of Poincar\'e
differential forms with algebraic integral curves of degrees $\mu$ in general position. For $d=3$ we calculate their dimension and show that their closures form components of $X_\infty$.

\begin{notation}
For the purpose of this section we consider

\begin{tabular}{ll}
$F \in K[x,y,z]$ & a homogeneous Polynomial of degree $e$\\
$C_F := \{ F=0 \}$ & the corresponding algebraic curve in $\PP^2$\\
\end{tabular}
\end{notation}

\begin{rem}
$C_F$ is an integral curve of $Pdx+Qdy$ if and only if $PF_y-QF_x=0$ on
$C_F$. This is the case if and only if there exists
a unique homogeneous polynomial $K_F \in K[x,y,z]$ of degree $d-1$ such that
\[
	PF_y-QF_x = -FK_F.
\]
We put the minus on the right side to simplify the formulas later on.
\end{rem}

\begin{defn}
In this situation $K_F$ is called the {\sl cofactor} of $F$.
\end{defn}

The importance of algebraic integral curves and their cofactors is given by

\begin{thm}[Darboux] \xlabel{tIntegrate}
Let $F_1,\dots,F_n \in K[x,y,z]$ be homogeneous Polynomials, such that 
the corresponding
algebraic curves  $C_1,\dots, C_n$ are integral curves of $Pdx+Qdy$  
with cofactors $K_1,\dots,K_n$. If there exist scalars $\alpha_i \in K$ such that
\[
	Q_x - P_y = \sum_{i=1}^n \alpha_i K_i
\]
then
\[
	\mu = \prod_{i=1}^n F_i^{\alpha_i}
\]
is an integrating factor for $Pdx + Qdy$. Such differential forms are called {\sl Darboux integrable}.
\end{thm}

\begin{proof} \cite{Darboux}
\end{proof}

If the algebraic integral curves have the right degrees, one doesn't even have to consider the cofactors:

\begin{thm}[Christopher, \.Zo\l\c{a}dek]\xlabel{tFamily}
Let $C_1,\dots,C_n \subset \PP^2$ be smooth algebraic curves of degrees $e_1, \dots, e_n$ with
\begin{enumerate}
\item  no curve tangential to the line at infinity, \xlabel{oneinfinity}
\item  no two curves tangential to each other, \xlabel{twotangent}
\item no two curves meeting the line at infinity at the same point, \xlabel{twoinfinity}
\item no three curves meeting in the same point. \xlabel{threemeet}
\end{enumerate}

Then the vector space of degree $d = e_1+\dots+e_n-1$ differential
forms $Pdx+Qdy$ that have these integral curves is $n$ dimensional. Further more all these
differential forms are Darboux integrable.
\end{thm}

\begin{proof} $Pdx+Qdy$ has $C_1,\dots,C_n$ as integral curves with cofactors $K_1,\dots,K_n$ if and only if
$(K_1,\dots,K_n,-Q,P) \cdot M = 0$ where
\[
	M = 
	\begin{pmatrix}
	F_1 & \hdots  & 0 \\
	\vdots & \ddots & \vdots \\
	0 & \hdots & F_n\\
	F_{1x} & \hdots & F_{nx} \\
	F_{1y} & \hdots & F_{ny} \\
	\end{pmatrix}
\]
and the $F_i$ are defining equations of the $C_i$. The conditions of the theorem are satisfied if and only
if $M$ has full rank everywhere on $\PP^2$.  In this case a presentation of the kernel of $M$
is given by the Buchsbaum-Rim-complex \cite[Apendix A2.6]{Ei95}. This complex gives the number of 
degree $d$ differential forms and writes them in terms of $n \times n$ minors of $M$. A short calculation then gives that $Q_x-P_y$ is in the span of $K_1,\dots,K_n$. For a proof involving only Hilbert's Nullstellensatz and no complexes see \cite[Theorem 1]{withCurves} and \cite[Theorem 1]{enoughCurves}.
\end{proof}

We now pass from general differential forms $Pdx+Qdy$ to Poincar\'e differential forms:

\begin{lem} \xlabel{lDepend}
Let $Pdx+Qdy$ a Poincar\'e differential form that is Darboux integrable with respect to 
the algebraic integral curves $C_1,\dots,C_n$. Let $F_1,\dots,F_n$ be the defining equations and $K_1,\dots,K_n$ the cofactors. If $a \in \PP^2$ is a point not lying on any $C_i$ satisfying $P(a)=Q(a)=0$ then $Q_x(a)=P_y(a)$.
\end{lem}

\begin{proof}
For each curve $C_i$ we have
\[
	P(a)F_{iy}(a)-Q(a)F_{ix}+F_i(a)K_i(a) = 0 \implies K_i(a) = 0
\]
since $P(a)=Q(a)=0$ and $F_i(a) \not=0$. Consequently
\[
	Q_x(a)-P_y(a) = \sum_{i=0}^n \alpha_i K_i(a) = 0
\]
since $Pdx+Qdy$ is Darboux integrable with respect to the $C_i$.
\end{proof}

\begin{defn}
We denote by $Y_\lambda$ with $\lambda = (\lambda_1,\dots,\lambda_n) \vdash d+1$ the
variety of differential forms that have $n$ integral curves of degree $\lambda_i$
satisfying the conditions of Christopher's Theorem. For the subvariety of Poincar\'e differential forms we use the notation $X_\lambda = Y_\lambda \cap \Pd$.
\end{defn}

\begin{prop}
Let $d=3$. Then 
\begin{enumerate}
\item $\overline{X_{(4)}}$ is a codimension $5$ component of $X_\infty$,  
\item $\overline{X_{(3,1)}}$ is a codimension $6$ component of $X_\infty$,  
\item $\overline{X_{(2,2)}}$, $\overline{X_{(2,1,1)}}$, $\overline{X_{(1,1,1,1)}}$  are codimension $7$ components of $X_\infty$.
\end{enumerate}
\end{prop}

\begin{proof}
The variety of degree $4$ smooth plane curves, that do not pass through the origin is $14$ dimensional. By Christopher's Theorem this variety contains an Zariski open subset, such that for each curve of this subset there is a $1$ dimensional family of differential forms $Pdx+Qdy$ that have this integral curve. For degree reasons the cofactor in this case has to be $0$. By Christopher's Theorem these differential forms are Darboux-integrable with $Q_x-P_y = 0$. This means that $Pdx+Qdy$ is exact and there exists a polynomial $H$ of degree $4$ with $P=H_x$ and
$Q=H_y$. Therefore each such differential form has a $1$-dimensional family of algebraic
degree $4$ integral curves.

Furthermore the
Poincar\'e differential forms $\Pdthree$ are defined by $6$ linear equations in the space of
all differential forms of degree $3$. By Lemma \ref{lDepend} at most $5$ of these are independent on the set of Darboux-integrable differential forms. 
Combining these dimension counts we obtain
\[
	\dim X_{(4)} \ge 14+1-1-5 = 9.
\]
Since $\dim \Pdthree = 14$ this means $\codim X_{(4)} \le 5$. On the other hand
\[
{\scriptstyle
( {x}-9 x^{2}+2 {x} {y}+2 y^{2}+2 {x} y^{2}-y^{3})dx+({y}+x^{2}+4 {x} {y}-6 y^{2}-2 x^{2} {y}-3 {x} y^{2})dy
}
\]
is a Hamiltonian Poincar\'e differential form with potential
\[
{\scriptstyle
     (1/2) x^{2}+(1/2) y^{2}-3 x^{3}+x^{2} {y}+2 {x} y^{2}-2 y^{3}+x^{2} y^{2}-{x} y^{3}.
}
\]
It has a $1$-dimensional family of smooth plane quartic integral
curves that do not pass through the origin and the tangent space of $X_{10}$ at this point
over $F_{23}$ is $5$. By semicontinuity arguments this proves that $\codim X_{(4)} \ge 5$ and that
$\overline{X_{(4)}}$ is a component of $X_\infty$.

The variety of smooth degree $1$ and degree $3$ curves have dimension $2$ and $9$ respectively. The variety of pairs that intersect transversally
is then $11$ dimensional. For an open subset of these pairs there is a $2$ dimensional set of 
differential forms that have these integral curves. Again by Lemma $\ref{lDepend}$ the codimension of the variety
of Poincar\'e differential forms in this set is at most $5$. Since 
\[
{\scriptstyle ( {x}-10 x^{2}+2 {x} {y}-2 y^{2}-8 x^{2} {y}-11 {x} y^{2}+3 y^{3})dx+({y}-5 x^{2}-4 {x} {y}-3 y^{2}-x^{3}-9 x^{2} {y}+10 {x} y^{2}-8 y^{3})dy}
\]
is a Poincar\'e differential form with exactly one degree $1$ and one degree $3$ integral curve 
\[
{\scriptstyle 10 y+{z}=0  \quad{\displaystyle \text{and}}\quad 10 x^{3}-11 x^{2} {y}+6 {x} y^{2}-9 y^{3}+10 x^{2} {z}+2 y^{2} {z}+2 {y} z^{2}+z^{3}=0}
\]
over $F_{23}$ satisfying the conditions of Christopher's Theorem, we have
\[
	\dim X_{(3,1)} \ge 11 + 2 - 0 -5 = 8.
\]
Since the tangent space of $X_{10}$ at this point has codimension $6$, this proves that $\overline{X_{(3,1)}}$ is a codim $6$ component of $X_\infty$.

For the remaining components we observe that the variety of degree $2$ plane curves is $5$ dimensional and that
\begin{align*}
	7 &= (5+5)+2-5 \\ &= (5+2+2)+3-5 \\ &= (2+2+2+2)+4-5.
\end{align*}

Now the Poincar\'e differential form
\[
{\scriptstyle
( {x}-2 x^{2}-11 {x} {y}+9 y^{2}+3 x^{3}+2 x^{2} {y}-11 {x} y^{2}+11 y^{3})dx+({y}+4 x^{2}-3 {x} {y}+y^{2}+2 x^{3}-11 x^{2} {y}+11 {x} y^{2}-7 y^{3})dy} 
\]
has two integral conics
\[
{\scriptstyle -6 x^{2}-7 {x} {y}-10 y^{2}+6 {x} {z}-11 {y} {z}+z^{2} =0 \quad{\displaystyle \text{and}}\quad 2 x^{2}+8 {x} {y}+11 y^{2}-4 {x} {z}-8 {y} {z}+z^{2}=0}
\]
in satisfying the conditions of Christopher's Theorem over $\FF_{23}$.
The Poincar\'e differential form
\[
{\scriptstyle
( {x}-3 x^{2}+11 {x} {y}+5 y^{2}+8 x^{3}+2 x^{2} {y}+8 {x} y^{2}-4 y^{3})dx+({y}-5 x^{2}-3 {x} {y}-11 y^{2}+7 x^{3}+5 x^{2} {y}-3 {x} y^{2})dy} 
\]
has two integral lines and one integral conic
\[
{\scriptstyle -6 x+{z} =0,\quad -3 x+3 y+{z}=0 \quad{\displaystyle \text{and}}\quad  x^{2}+5 {x} {y}-2 y^{2}-4 {x} {z}-5 {y} {z}+z^{2} =0}
\]
in satisfying the conditions of Christopher's Theorem over $\FF_{23}$.
The Poincar\'e differential form
\[
{\scriptstyle
( {x}-11 x^{2}+3 {x} {y}-2 y^{2}-3 x^{3}+9 x^{2} {y}-5 {x} y^{2}-11 y^{3})dx+({y}+5 x^{2}+8 {x} {y}-11 y^{2}+4 x^{3}-4 x^{2} {y}+8 {x} y^{2}+7 y^{3})dy
}
\]
has four integral lines 
\[ 
{\scriptstyle 3 x+8 y+{z}=0, \quad -2 x+9 y+{z}=0, \quad {x}+5 y+{z}=0 \quad{\displaystyle \text{and}}\quad 10 x-6 y+{z}=0}
\]
in satisfying the conditions of Christopher's Theorem over $\FF_{23}$. 

Moreover the tangent spaces of $X_{10}$
in these points over $\FF_{23}$ are all of codimension $7$.
\end{proof}

\begin{rem}
All examples in this section were found with the program \cite{strudelweb}.
\end{rem}

\section{Rationally Reversible Systems} \xlabel{sReversible}
\nosubsections

A second type of centers has been considered by \zoladek:

\begin{defn}
A differential form $Pdx+Qdy$ is called rationally reversible, if there exist a rational map
\[
	\Phi \colon \CC^2 \to \CC^2
\]
and a second differential form $P'dx + Q'dy$ such that 
\[
	\Phi^*(P'dx+Q'dy) = \mu (Pdx + Qdy)
\]
with $\mu$ a suitable polynomial.
\end{defn}

This definition is useful because of

\begin{prop}[\zoladek]
Let $\mu (Pdx + Qdy) = \Phi^*(P'dx+Q'dy) $ be a rationally reversible differential form over $\RR$
and $O \in \RR^2$ a ramification point of $\Phi$. Furthermore doenote by $\Gamma_\Phi$ the ramification curve of $\Phi$ and $\Phi(\Gamma_\Phi)$ the branch curve. If
\begin{enumerate}
\item there exists a neighborhood $U \subset \RR^2$ of $O$ such that the boundary of $\Phi(U)$
contains part of the branch curve $\Phi(\Gamma_\Phi)$,
\item $P'dx + Q'dy$ does not vanish at $\Phi(O)$ and
\item the integral curve of $P'dx+Q'dy$ through $\Phi(O)$ is tangent to the branch curve
$\Phi(\Gamma_\Phi)$ from 
the outside of $\Phi(U)$,
\end{enumerate}
then $O$ is a center of $Pdx+Qdy$.
\end{prop}

\begin{proof}(\zoladek)
In this case the integral curves of $Pdx+Qdy$ close to $O$ are preimages of compact curves in
$\Phi(U)$.
\end{proof}

\begin{defn}
In the situation above, let $\Gamma_{fold} \subset \Gamma_\Phi$ be the union of those components
of $\Gamma_\Phi$ that contain $O$. $\Gamma_{fold}$ is then called the {\sl fold curve} of $\Phi$.
\end{defn}

The first examples of rationally reversible differential forms are those which are symmetric with respect to the
$x$ axis. These were already considered by Poincar\'e. In this case $\Phi = (x,y^2)$. In \zoladeks list
this ist the family $CR_{1}$.

\newcommand{\xxaxisR}{X_{x-axis,\RR}}
\newcommand{\xxaxisC}{X_{x-axis,\CC}}
\newcommand{\xxaxis}{X_{x-axis}}
\newcommand{\xsym}{X_{sym}}

We now show, that the Poincar\'e differential forms which are symmetric to a line though the
origin, form a codimension $6$ component $\xsym$ of $X_\infty$.

\begin{prop}
Let $\xxaxisC \subset \Pd$ be the variety of complex valued Poincar\'e differential forms that are mirror symmetric with respect to the $x$-axis,
i.e. those $Pdx+Qdy$ that satisfy $P(x,-y,z)=P(x,y,z)$ and  $Q(x,-y,z)=-Q(x,y,z)$. Then $\xxaxis \subset X_\infty$.
\end{prop}

\begin{proof}
$\xxaxisC$ is a vector subspace of $\Pd$, since 
symmetry with respect to the $x$ axis is a linear condition. For real symmetric differential forms $\xxaxisR \subset \xxaxisC$ the focal values $s_j$ vanish. (see for example \cite[p. 42, Satz 5.2.1]{Prell}). Since the $s_j$ are holomorphic on $\xxaxisC$ they have to vanish on all of $\xxaxisC$.
\end{proof}


\begin{rem}
The same result can also be obtained
analytically  by introducing complex polar coordinates \cite{vWahlGap}. 
\end{rem}

\begin{rem}
Over the real field mirror symmetry implies stability even for differential equations of the form
\[
	y' = - \frac{x^{2n-1} + p(x,y)}{y^{2n-1}+q(x,y)}
\]
 \cite[p. 42, Satz 5.2.1]{Prell}.\end{rem}

\begin{prop}
Let $\xsym \subset X_\infty \subset \Pdthree$ be the variety of mirror symmetric Poincar\'e differential forms
of degree $3$. Then $\xsym$ is a component of $X_\infty$ and $\codim \xsym = 6$. 
\end{prop}

\begin{proof}
$\xxaxis$ is defined by the
linear conditions $p_{ij} = 0$ for $j$ odd and $q_{ij}=0$ for $j$ even. 
Therefore $\codim \xxaxis = 7$. 

By rotating the symmetry axis about
the origin, i.e. by applying an element of $\SO(2,\CC)$, we obtain a $1$-parameter family of mirror symmetric differential forms in $X_\infty$ for each element $Pdx+Qdy \in \xxaxis$. Since for a general $Pdx+Qdy \in \xxaxis$ this family is not constant, 
this shows $\codim \xsym \le 6$.
 
Now 
 \[
 ( {x}+11 x^{2}+y^{2}{-9 {x} y^{2}})dx+({y}{-6 {x} {y}}+x^{2} {y}-5 y^{3})dy
 \]
is mirror symmetric. In characteristic $23$ it can be checked with Frommer's algorithm that
$\codim T_{x,X_{10}} = 6$. Since $\xsym \subset X_{10}$, this shows $\codim \xsym \ge 6$. 
It also shows that the component of $X_\infty$ that passes through $x$ has codimension
at least $6$. Since $\xsym$ is irreducible of codimension at most $6$ and contains $x$, it 
must be this component.  

\end{proof}

Next we consider Zoladek's family $CR_{11}$.

\begin{defn}
Let $\Phi = (A^2/B,A/C)$ with $A$,$B$ and $C$ linear polynomials in $x$ and $y$, and 
\[
	P'dx+Q'dy =  -2x(ky+lx+mxy)dx +  
              y(ky+nx+py^2+qxy+rxy^2) dy.
\]
Then
\[
	Y_{CR11} = \{ Pdx+Qdy \suchthat \exists A,B,C,k,l,m,n,p,q,r,\mu \colon
	\mu(Pdx+Qdy) = \Phi^*(P'dx+Q'dy) \}
\]
is called \zoladeks family $CR_{11}$. We also set $X_{CR11} := Y_{CR11} \cap \Pdthree$. 
\end{defn}

\begin{prop}
$\overline{X_{CR11}}$ is a codim $7$ component of $X_\infty$.
\end{prop}

\begin{proof}
Since $Y_{CR11}$ is rational and the focal values are
holomorphic functions on the open set of differential forms with an elementray center that vanish on
the real points of $Y_{CR11}$ by \zoladeks Theorem $1$ in \cite{ZoladekRational} they must also vanish on all of $Y_{CR11}$. 

Also  $Y_{CR11}$ is a codimension $7$ family of differential forms by the same theorem. Therefore $X_{CR11}$ has codimension at most
$7$. Now
\[
{\scriptstyle
(-8 x^{3}+12 x^{2} {y}+6 {x} y^{2}-18 y^{3}+x^{2}-2 {x} {y}-6 y^{2}+{x})dx + 
(5 x^{3}+3 x^{2} {y}-11 {x} y^{2}+11 y^{3}+2 x^{2}-15 {x} {y}-15 y^{2}+{y})dy
}
\]
is an element of $X_{CR11}$ with a codim $7$ tangent space in $X_{17}$ over $\FF_{37}$.
This proves the proposition.
\end{proof}

\begin{rem}
For a random subset of the codim $7$ points we found in our experiments we also checked wether they lie on $X_{CR11}$. This was the case for about one quarter of the points, which is consistent with our conjecture that there are exactly $4$ codimension $7$ components of which $\overline{X_{CR11}}$ is one.
\end{rem}

\section{\zoladeks List} \xlabel{sList}
\nosubsections

In an amazing work \zoladek classified all rationally reversible cubic systems that are not Darboux integrable (\cite{ZoladekRational}, \cite{ZoladekCorrection}).  A list of the $17$ such centers is presented in \cite{ZoladekRational} and some minor mistakes are corrected in \cite{ZoladekCorrection}. In this section we compare these results to our experiments.

\begin{prop}
Let $Y_{CRi} \subset \Vdthree$ be the family $CR_i$ from \zoladeks list, and $X_{CRi} = Y_{CRi} \cap \Pdthree$
the subspace of Poincar\'e differential forms in this family. Then
\begin{itemize}
\item $\codim X_{CRi} \ge 6$ for $i=1$,
\item $\codim X_{CRi} \ge 7$ for $i = 5,7,11,12,16$,
\item $\codim X_{CRi} \ge 8$ for $i = 2,4,6,8,13,14$,
\item $\codim X_{CRi} \ge 9$ for $i = 3,9,10,15$,
\item $\codim X_{CRi} \ge 10$ for $i = 17$.
\end{itemize}
\end{prop}

\begin{proof}
Points on $X_{CRi}$ that have the have tangent spaces of the above codimension can be found by the following method. Choose a random differential Form in $x \in Y_{CRi}$ over a finite field $\Fp$. If there is a $\Fp$-rational point on the fold curve, where the corresponding differential form degenerates, apply a change of coordinates $\phi$ that moves this point to the origin. If the linear part of $x' = \phi(x)$ admits a $\Fp$-rational coordinate change $\psi$ such that $\psi(x')$ is in $\Pdthree$ set $y = \psi(x')$ and
calculate the tangent space with Frommer's algorithm. This methods usually gives an example after $4$ or $5$ trials. The proposition then follows from semi-continuity.
\end{proof}

\begin{rem}
In \zoladeks list of \cite{ZoladekCorrection} there are some minor misprints. For $CR_3$ one has to replace $(2a^2-b)$ in the formula for $F$ by $(2a-b^2)$. For $CR_5$ the fromulas for $\dot{x}$ and $\dot{y}$ have to be exchanged and $-l$ has to be replaced by $-ly$. For $CR_{12}$ the term $-pqT$ must be replaced by $-yqT$. In every case either the list in \cite{ZoladekRational}, the one in \cite{ZoladekCorrection} or a combination of the two is correct.
\end{rem}

\begin{rem}
Because of the new formula for $\eta$ given in \cite{ZoladekCorrection} the fold curve in this case
changes to $\Gamma= \{x-2ay+1-a=0\}$. 
\end{rem}

\begin{rem}
Because $\codim(X_{CRi}\subset \Pdthree) \ge \codim(Y_{CRi} \subset \Vdthree)$ our examples for the family $CR_{16}$
show that the codimension $5$ given in \zoladeks list can not be correct.
\end{rem}

We will now take a closer look at those families that could be of codimension $7$ in the
space of Poincar\'e differential forms. As in \zoladeks papers we set
\[
	T = x+y+c \quad\quad\text{and} \quad\quad T_2= ax^2+bxy+cy^2+dx+ey+1.
\]
We have found some previously unknown conic integeral curves in some of \zoladeks rationally reversible families, proving that they are also Darboux integrable:


\begin{prop} \label{pCR16}
A general differential form of \zoladeks family $CR_{16}$
\begin{align*}
{\scriptstyle
Pdx+Qdy =} 
 & {\scriptstyle
 -\bigl(-y (k y+l x)-(p y+q x)*T_2-(2 a x+b y+d) (-q x^2+(n-p) x y+m y^2)\bigr) dx
 }\\
 &{\scriptstyle
+\bigl(-x (k y+l x)-(m y+n x) T_2+(b x+2 c y+e) (-q x^2+(n-p) x y+m y^2)\bigr) dy
}
\end{align*}
has a conic integral curve
\[
	\bigl\{ (mq-np)T_2+(kq- lp)x-(kn-lm)y = 0 \bigr\}
\]
and is Darboux integrable with $X_{CR16} \subset \overline{X_{(2,1,1)}}$.
\end{prop}

\begin{proof}
That the conic above is an integral curve can be checked by a straight forward calculation. By \zoladeks Remark
 $1$ in \cite{ZoladekRational} the differential forms of this system also have 
 two integral lines given by the equation $qx^2+(p-n)xy-my^2=0$. By substituting random values of $\Fp$ for the parameters of this family one easily obtains examples that satisfy the conditions Christophers Theorem \ref{tFamily} over $\Fp$. The proposition then follows by semi-continuity.
\end{proof}


\begin{prop}
A general differential form of \zoladeks family $CR_{12}$
\begin{align*}
{\scriptstyle
Pdx+Qdy =} 
 & {\scriptstyle
 -\bigl(-k y^2+2 (q-m) x y+2 r x T-2 l y^3-y q T-r T^2-2 n y^2 T-2 p y T^2\bigr) dx
 }\\
 &{\scriptstyle
     + 2 \bigl((m-k-q) x y-(m+r) x T+l y^3+(p-n) y T^2+(n-l) y^2 T-p T^3\bigr) dy
}
\end{align*}
has a conic integral curve
\[
\{a_{20} x^2+a_{11} xy + a_{02}y^2 + a_{10}x + a_{01}y + a_{00}=0\}
\]
where the $a_{ij}$ satisfy linear system of the equations
\[
\begin{pmatrix}0&
     {-c}&
     0&
     0&
     1&
     0\\
     -2 {p} {c}-m-r&
     1/2 r&
     0&
     {p}&
     0&
     0\\
     -2 {n} {c}-k+2 m-q+{r}&
     -m+1/2 q-r&
     {r}&
     {n}&
     0&
     0\\
     -2 {l} {c}+{k}-m+{q}&
     -1/2 k+{m}-q&
     -m+{q}&
     {l}&
     0&
     0\\
     {-c^{2}}&
     0&
     0&
     0&
     0&
     1\\
     \end{pmatrix}
     \begin{pmatrix} a_{20} \\ a_{11} \\ a_{02} \\ a_{10} \\ a_{01} \\ a_{00} \end{pmatrix}
     = 0.
\]
Furthermore $Pdx+Qdy$ is Darboux integrable and $X_{CR12} \subset \overline{X_{(2,1,1)}}$.
\end{prop}
 
\begin{proof}
Same as for Proposition \ref{pCR16}, with \zoladeks integral lines given by the equation $ky^2+qTy+rT^2=0$. (In \zoladeks paper the equation is $ky^2+pTy+qT^2=0$, but this is a misprint as can be checked by a straight forward computation).
\end{proof}

For the families $CR_7$ and $CR_5$ the integral curves we found do not satisfy the conditions of Christophers Theorem. To prove that they are subfamilies of $\overline{X_{(2,2)}}$ we need the following technical lemma:


\begin{lem} \label{2x3}
Let $Pdx+Qdy$ be a differential form of degree $3$, 
\[
	M = \begin{pmatrix} F_x & G_x \\ F_y & G_y \\ aF & bG \end{pmatrix}
\]
a matrix with $F$, $G$ quadratic polynomials and $a,b \in \CC \backslash \{0\}$ such that
\begin{enumerate}
\item $(-Q,P,Q_x-P_y) \cdot M = 0$ and 
\item $M$ drops rank in only finitely many points of $\PP^2_\CC$.
\end{enumerate}
Then $Pdx+Qdy$ is Darboux integrable and an element of $\overline{X_{(2,2)}}$.
\end{lem}

\begin{proof}
In this situation $F$ and $G$ are integral curves of $Pdx+Qdy$ with cofactors
$a(Q_x-P_y)$ and $b(Q_x-P_y)$. In particular $Pdx+Qdy$ is Darboux integrable.

If $M$ drops rank in codimension $2$ then the syzygies of $M$ are generated by
\[
	s = \left( 
		\begin{vmatrix} F_y & G_y \\ aF & bG \end{vmatrix},
		-\begin{vmatrix} F_x & G_x \\ aF & bG \end{vmatrix},
		\begin{vmatrix} F_x & G_y \\ F_y & G_y\end{vmatrix}
	   \right)
\]
For degree reasons the syzygy $(-Q,P,Q_x-P_y)$ has to be a scalar multiple of $s$, in particular
\[
	Q = -c \begin{vmatrix} F_y & G_y \\ aF & bG \end{vmatrix}
	\quad \quad \text{and} \quad \quad
	P = -c \begin{vmatrix} F_x & G_x \\ aF & bG \end{vmatrix}.
\]
We now deform $F$ and $G$ to $\tilde{F}$ and $\tilde{G}$ such that $C_{\tilde{F}}$ and $C_{\tilde{G}}$
satisfy the conditions of Christopher's Theorem. This is possible, since these conditions are Zariski open
on the set of all pairs of conics. Then $C_{\tilde{F}}$ and $C_{\tilde{G}}$ are integral curves of $\tilde{P}dx+\tilde{Q}dy$ with
\[
	\tilde{Q} = -c \begin{vmatrix} \tilde{F}_y & \tilde{G}_y \\ a\tilde{F} & b\tilde{G} \end{vmatrix}
	\quad \quad \text{and} \quad \quad
	\tilde{P} = -c \begin{vmatrix} \tilde{F}_x & \tilde{G}_x \\ a\tilde{F} & b\tilde{G} \end{vmatrix}
\]
and cofactors
\[ 
ca\begin{vmatrix} \tilde{F}_x & \tilde{G}_y \\ \tilde{F}_y & \tilde{G}_y\end{vmatrix}
\quad\quad\text{and}\quad\quad
cb\begin{vmatrix} \tilde{F}_x & \tilde{G}_y \\ \tilde{F}_y & \tilde{G}_y\end{vmatrix}.
\]
By Christophers Theorem $\tilde{P}dx+\tilde{Q}dy \in X_{(2,2)}$ is Darboux integrable. Now $\tilde{P}dx+\tilde{Q}dy \in X_{(2,2)}$ deforms to
$Pdx+Qdy$, which is therefore in $\overline{X_{(2,2)}}$. \end{proof}


\begin{prop} \xlabel{pCR7}
A general differential form of \zoladeks family $CR_7$
\[
{\scriptstyle
Pdx+Qdy = 
     -\bigl(n x+k y+n T+(m-l) x^2 y+p x^2 T+m x y T+p x T^2\bigr) dx +
     x \bigl(-(n+k)+(l-m) x y-(l+p) x T\bigr) dy
}
\]
has a conic integral curve
\[
 C_F = \bigl\{(lp-mp)x^2 + (lm-m^2+lp-mp)xy + c(lp-mp)x + (-mn+kp)=0\bigr\}.
 \]
 and is Darboux integrable with $X_{CR7} \subset \overline{X_{(2,2)}}$.
 \end{prop}
 
 \begin{proof}
 That the conic $C_F$ above is an integral curve can be checked by a straight forward calculation. 
 By \zoladeks Remark $1$ in \cite{ZoladekRational} the differential forms of this system also have another integral conic $C_G$ with $G=k+lTx$. With $a = \frac{l-m}{2l-m}$ and $b= \frac{l}{2l-m}$
 we obtain a matrix $M$ as in Lemma \ref{2x3}. Condition (1) of Lemma \ref{2x3} follows from a direct computation and condition (2) can be verified for random values of the coefficients over a finite field, since the rank condition is Zariski open on the variety of all such matrices. The claim then follows from Lemma \ref{2x3}.
 \end{proof}


The most involved argument for Darboux integrability is needed for the family $CR_5$:

\begin{prop}
A general differential form of \zoladeks family $CR_5$
\[
{\scriptstyle
Pdx+Qdy = 
     \bigl(-k x y^2-l y+m x^2 y-(n x y+p+q T x) (2 x+y+c)\bigr) dx 
     -x (l+p+m c x+(k+n) x y+m x^2+q T x) dy
}
\]
has two conic integral curves
\[
	\{a_{20} x^2+a_{11} xy + ca_{20}x + a_{00}=0\}
\] 
where $a_{ij}$ are such that the matrix
\[
M =
\begin{pmatrix}{a_{20}}&
      {q}&
      -m-q\\
      {-a_{11}}&
      {m}-n-q&
      {k}+{n}+{q}\\
      {a_{00}}&
      {p}&
      -l-p\\
      0&
      {-a_{20}}&
      {a_{11}}\\
      \end{pmatrix}
\]
has rank at most $2$. $Pdx+Qdy$ is
Darboux integrable with respect to these quadrics and $X_{CR5} \subset \overline{X_{(2,2)}}$.
\end{prop}

\begin{proof}
Indeed for a general member of \zoladeks family $CR_5$ we may after row and column
operations assume that $M$ has the form
\[
\begin{pmatrix}
 0 & 1 & 0 \\
 0 & 0 & 1 \\
 L & 0 & 0 \\
 Q & 0 & 0
\end{pmatrix}
\]
with $L$ linear in $a_{ij}$ and $Q$ quadratic in $a_{ij}$. So $M$ drops rank in two points and
we have two quadrics. With a computer algebra program one can check that the equations given define an ideal, that contains the ideal that describes the set of all quadric integral curves. Consequently the quadric of the proposition are integral.

Now consider the matrix
\[
{\scriptsize
 M = \begin{pmatrix}2 {y} {m}-{y} {n}&
      2 {m} {t} {c}-{n} {t} {c}+4 {x} {m}-2 {x} {n}\\
      2 {x} {m}-{x} {n}&
      0\\
      -{x} {q} {t} {c}+{x} {y} {m}-{x} {y} {n}-x^{2} {q}-{x} {y} {q}-{p} t^{2}&
      {x} {m} {t} {c}+{x} {q} {t} {c}+{x} {y} {k}+x^{2} {m}+{x} {y} {n}+x^{2} {q}+{x} {y} {q}+{l} t^{2}+{p} t^{2}\\
      \end{pmatrix}
}
\]
which we obtained as syzygy matrix of $(-Q,P,Q_x-P_y)$. It does indeed satisfy $(-Q,P,Q_x-P_y) \cdot M = 0$. To show that it can be transformed into the form of Lemma \ref{2x3}, we multiply the last line
by $\alpha$ and subtract its derivatives from the first two lines. Looking at the coefficients
of the resulting linear entries we obtain a $2 \time 6$-matrix
\[
{\scriptsize
M' = \left(
\begin{array}{c|c}
2 q&
      4 {m} {\alpha}-2 {n} {\alpha}-2 m-2 q\\
      2 {m} {\alpha}-{n} {\alpha}-m+{n}+{q}&
      -k-n-q\\
      2 {m} {\alpha}-{n} {\alpha}-m+{n}+{q}&
      -k-n-q\\ \hline
      0&
      0\\
      {q} {c}&
      2 {m} {\alpha} {c}-{n} {\alpha} {c}-{m} {c}-{q} {c}\\
      0&
      0\\
\end{array}
\right).
}
\]
This matrix has a kernel if and only if the first $2\times 2$ minor 
\[
	\phi =
      (8 m^{2}-8 {m} {n}+2 n^{2})(\alpha-\alpha^2)
      -2 m^{2}+2 {m} {n}-2 {k} {q}
\]
vanishes, since all other rows are dependent on the first two. For generic choices of $k$, $m$, $n$ and $q$ the polynomial $\phi$ has two zeros $\alpha_1 \not= \alpha_2 \not= 0$. For each of these values $M'$ has a kernel $\beta_i = (\beta_{i1},\beta_{i2})^T$. By construction we obtain
\[
	M \cdot \begin{pmatrix} \beta_{11} & \beta_{21} \\ \beta_{12} & \beta_{22} \end{pmatrix}
	=  \begin{pmatrix} F_x & G_x \\ F_y & G_y \\ \frac{1}{\alpha_1}F & \frac{1}{\alpha_2}G \end{pmatrix}
\]
with appropriate quadric polynomials $F$ and $G$. 

By substituting random values for $c$, $k$, $m$, $n$ and $q$ over a finite field $\FF_p$ we easily find an example where  $\alpha_1$ and $\alpha_2$ are $\FF_p$-rational and the constructed matrix 
satisfies the rank condition of Lemma \ref{2x3}. The proposition then follows by semi continuity.
\end{proof}

\begin{cor}
The only family of rationally reversible cubic centers of codimension $7$ in $\Pdthree$ whose general member is not Darboux integrable with integral curves of degree at most two is $CR_{11}$.
\end{cor}

\begin{proof}
By looking at random points of $CR_{11}$ over $\Fp$ one can easily find an example that is
not Darboux integrable with integral curves of degree at most two. $CR_1$ is the codim $6$ component of symmetric differential forms. All other families of
\zoladeks list are either Darboux integrable or have codimension at least $8$.
\end{proof}

\begin{rem}
It would be interesting to find a geometric explanation for the extra conic integral curves
we have found.
\end{rem}

\section{Necessary equations} \xlabel{sNecessary}
\nosubsections

As a last application of our methods we estimate the number of focal polynomials needed to
define $X_\infty$.

\begin{prop}
Let $x \in \AZ^n$ be a point with $s_1(x) = \dots = s_k(x) = 0 \mod p$, $s_{k+1}(x) \not=0$ and
$\codim T_{X_k,x} = k$ where $X_k = V(s_1,\dots,s_k)$.Then there exists a component of $X_{k}$
that is not contained in $X_{k+1}$ in characteristic $0$.
\end{prop}

\begin{proof}
Since the maximal codimension of $X_k$ is $k$, $x$ lies on a component $X'$ of $X_k$ 
in characteristic $p$ that is smooth in a neighborhood of $x$ and of codimension $k$. By a theorem
of Schreyer (\cite{smallFields}, \cite{newFamily}), this component lifts to a component $\tilde{X}'$ over an algebraic number field. 
Now $x$ does not lie in $X_{k+1}$ considered as a scheme over $\ZZ$. Therefore $s_k$ can not
vanish on all of $\tilde{X}' \subset X_{k}$.
\end{proof}

\begin{example}
For
\[
{\scriptstyle
  Pdx + Qdy = (x-4x^2+7xy+11y^2+7x^3+7x^2y-6xy^2-3y^3)dx
                       +(y-3x^2-5xy-4y^2-10x^3-4x^2y+11xy^2-6y^3)dy
}
\]
we have $s_1(P,Q) = \dots = s_9(P,Q) = 0 \mod 23$ and $s_{10}(P,Q) \not= 0 \mod 23$. Therefore
$m(3) \ge 10$. This example was found with the program \cite{strudelweb}. Notice that currently the best known bound is $m(3) \ge 11$ which requieres a more ingeneous proof \cite{ZoladekEleven}, \cite{ChristopherEleven}.
\end{example}

For upper bounds we have to fall back on heuristic methods:

\begin{heu}
If $X_k \not= X_{k+1}$ there exists a component $X' \subset X_k$ that is not completely
contained in $X_{k+1}$. Since all components of $X_k$ have codimension at most $k$ one would expect to find at least one point of $X'$ with probability
\[
	1 - \left( 1 - \frac{1}{p^k} \right)^N.
\]
If this probability is large and we do not find such a point after $N$ trials we conjecture $X_k = X_{k+1}$.
\end{heu}

\begin{example}
In Experiment \ref{eHom2} we expect to find at least one point on $X_3-X_4$ with probability
\[
	1 -  \left( 1 - \frac{1}{11^3} \right)^{100,000} \approx 99,945\%.
\]
if such points exist. Since we didn't find any points on $X_3-X_4$ this leads us to estimate $m(2) = 3$. Indeed this is the correct value.
\end{example}

\begin{example}
In Experiment \ref{eHom3} we would expect to find at least one point on $X_5-X_6$ with probability
\[
	1 -  \left( 1 - \frac{1}{17^5} \right)^{68,000,000} \approx (100- 1.5 \times 10^{-19}) \%.
\]
if such points exist. Since we did not find any points on $X_5 - X_6$ we conjecture that $m(3_{\text{homogeneous}}) = 5$. Indeed the intersection of condition I and condition II in \cite{LSHom3} is generated by the first $5$ focal polynomials.
\end{example}

\begin{example}
To obtain heuristic evidence for $m(3) = 11$ one would need to evaluate
\[
	N = \frac{\ln(0.05)}{\ln(1 - \frac{1}{29^{11}})}
	   \approx -\ln(0.05)29^{11} 
	\approx 3.7 \times 10^{16}
\]
random points. Unfortunately this about $10000$ times more than we can currently manage.
\end{example}


\begin{appendix}

\section{Frommer's algorithm}
\nosubsections
\begin{center}{\sc Hans-Christian Graf v. Bothmer and Martin Cremer} \end{center}

In this appendix we state Frommer's algorithm and show how one can use it to define the focal polynomials $\delta_{j}s_j$, evaluate them over finite fields and
calculate tangent spaces to their vanishing sets without ever writing them down explicitly.

\begin{alg}[Frommer, Moritzen] \xlabel{aFrommersAlgorithm}
Let $Pdx+Qdy$ be a Poincar\'e differential form on $\PP^2_K$. 
For $n \le 2k+2$ calculate $c_{l,n-l}$
for $0\le l \le n$ and $a_{n-l,l}$ according to the following formulas
\[
	c_{l,n-l} := \sum_{2 \le i+j \le n} \bigl(- (n-l-j+1)p_{ij}a_{l-i,n-l-j+1} + (l-i+1)q_{ij}a_{l-i+1,n-l-j} \bigr)
\]
\[
a_{n-l,l} := \left\{
	\begin{array}{cl}
		{\displaystyle -\sum_{i=1}^{\frac{n-l}{2}} \frac{c_{2i-1,n-(2i-1)}}{n-l}
		\prod_{j=i}^{\frac{n-l-2}{2}} \frac{n-2j}{2j}} & \text{for $l,n$ even, $l<n$} \\
		0 & \text{for $l=n$ even} \\
		{\displaystyle \frac{\alpha(\min(\frac{l-1}{2},\frac{n-l-1}{2}),n)}{2 \min(l,n-l)} 
		\left( \sum_{i=0}^{\frac{l-1}{2}} \frac{c_{n-2i,2i}}{\alpha(i,n)} - 
			\sum_{j=\frac{l-1}{2}+1}^{\frac{n}{2}} \frac{c_{n-2j,2j}}{\alpha(j,n)} \right)}
		& \text{for $l$ odd and $n$ even}\\
		{\displaystyle \sum_{i=0}^{\frac{l-1}{2}} \frac{c_{n-2i,2i}}{l} 
		\prod_{j=i}^{\frac{l-3}{2}} \frac{n-2j-1}{2j+1} }& \text{for $l,n$ odd} \\
		{ \displaystyle -\sum_{i=0}^{\frac{n-l-1}{2}} \frac{c_{2i,n-2i}}{n-l}
		\prod_{j=i}^{\frac{n-l-3}{2}} \frac{n-2j-1}{2j+1}} & \text{for $l$ even and $n$ odd} \\
	\end{array}
	\right.
\]
where the constants $\alpha(i,n)$ are defined as
\[
	\alpha(i,n) := \prod_{k=1}^{i} \frac{n-(2k-1)}{(2k-1)}.
\]
As start values use $a_{2,0}=a_{0,2}=1, a_{0,0}=a_{0,1}=a_{1,0}=a_{1,1}=0$. Also set $a_{i,j}=c_{i,j}=0$ if either $i<0$ or $j<0$. Then 
\[
	s_k(P,Q) = \frac{1}{2} \sum_{i=0}^{k+1} \frac{c_{2k+2-2i,2i}}{\alpha(i,2k+2)}.
\]
is called the {\sl $k$-th focal value of $Pdx+Qdy$}.
\end{alg}

\begin{prop}
If $\Char K = 0$ Frommer's algorithm is well defined and the formal power series $F(x,y) = \sum_{ij} a_{ij}x^iy^j \in K[[x,y]]$ satisfies 
\[
	\det \begin{pmatrix} F_x(x,y) & F_y(x,y) \\ P(x,y,1) & Q(x,y,1) \end{pmatrix} = \sum_{j=1}^\infty s_j(x^{2j+2}+y^{2j+2})
\]
\end{prop}

\begin{proof}
In characteristic zero all denominators that appear in Frommer's algorithm are invertible. Furthermore
the formula for $c_{l,n-l}$ involves only $a_{ij}$ with $i+j < n$ and the formula for
$a_{n-l,l}$ involves only $c_{ij}$ with $i+j \le n$. Therefore $s_k(P,Q)$ is
well defined. Moritzen shows in \cite{Moritzen} 
that $F$ satisfies the above formula. 

The idea and the first few steps of this algorithm appear in \cite{Frommer}. The explicit formulas above were found by Moritzen \cite{Moritzen}. For our implementation of Frommer's algorithm we followed \cite{Hoehn}. We used a random number generator from \cite{numericalRecipesC}. The source files of our program can be obtained from \cite{strudelweb}. 
\end{proof}

\begin{cor} \xlabel{cPolynomial}
The function $s_k \colon \Vd \to K$ is polynomial in $p_{ij}$ and $q_{ij}$ with rational coefficients,
i.e $s_k \in \QQ[p_{ij},q_{ij}]$.
\end{cor}

\begin{proof}
In Frommer's algorithm all formulas are algebraic and only integral numbers occur as denominators.
\end{proof}

\begin{cor}
If $\Char K=p$ then $s_k$ is well defined if $2k+2 <  p$.
\end{cor}

\begin{proof}
The denominators in the formulas for $c_{l,n-l}$ and $a_{n-l,l}$ and $s_k$ in Frommer's algorithms are products of natural numbers less or equal to $n \le 2k+2$.
\end{proof}

\begin{cor} \xlabel{cFiniteField}
Let $\Fp$ be a finite field of characteristic $p > 2k+2$ and $Pdx+Qdy$ a differential form over $\ZZ$. If $s_{k,\Fp}(\bar{P},\bar{Q})$ is the result of Frommer's algorithm over $\Fp$ and $s_{k,\QQ}(P.Q)$ the result of Frommer's algorithm over $\QQ$, then
\[
	s_{k,\Fp}(\bar{P},\bar{Q}) = 0 \iff \overline{\delta_k s_{k,\QQ}(P,Q)} = 0 
\]
where $\delta_k$ the smallest common denominator of the focal polynomial $s_{k,\QQ} \in \QQ[p_{ij},q_{ij}]$
and $\bar{P}dx+\bar{Q}dy$ is the reduction of $Pdx+Qdy$  to $\Fp$.
\end{cor}

\begin{proof}
By the argument above, $\delta_k$ is not divisible by $p$. 
\end{proof}

\begin{rem}
This allows us to find points on $X_i$ over a finite field with Frommer's algorithm without knowing the explicit polynomials $\delta_k s_k$.
\end{rem}

\begin{rem} \xlabel{rEpsilon}
For a Poincar\'e differential form $Pdx + Qdy$ over a field $K$ and  vector $(P',Q') \in \Pd$ 
we have
\[
	s_k(P+\epsilon P',Q+\epsilon Q') = s_k(P,Q) + \epsilon s_k'(P,Q) \in K[\epsilon]/(\epsilon^2)
\]
where $s_k'$ is the formal derivative of $s_k$ in the direction $(P',Q')$. Since Frommer's
algorithm works over any ring of characteristic $p \ge 2k+2$ this
allows us to calculate the tangent space of $X_k$ 
in the point $Pdx + Qdy$ without knowing the explicit polynomials $\delta_k s_k$.
\end{rem}

\begin{rem}
To speed up the search for rational points on $X_k$ over $\FF_p$,
we used the first focal polynomial (see Example \ref{eFirstFocal}) to calculate $q_{30}$ from the other values. This effectively
multiplies the number of points checked by $p$. The numbers reported in the main part of this paper are effective numbers, not actual numbers. Alternatively we could have used the actual number of points and calculated the codimensions in the hypersurface $s_1=0$ and then added $1$ to obtain the codimension in $\Pdthree$. The results are the same.
\end{rem}
 
\end{appendix}

\def\cprime{$'$}

\end{document}